\documentclass[14pt, letterpaper]{amsart}
\usepackage{amsmath,amsfonts,amssymb,amsthm}

\usepackage{graphicx,tikz}
\usepackage[all]{xy}

\usetikzlibrary{decorations.pathreplacing}
\usetikzlibrary{arrows}

\newtheorem{thm}{Theorem}[section]

\theoremstyle{definition}

\theoremstyle{remark}

\makeatletter
\let\c@equation\c@thm
\makeatother
\numberwithin{equation}{section}

\begin{document}

\title{On some linear combination of two contiguous hypergeometric functions.
}

\author{Imane GHANIMI
}

\dedicatory{Department of Mathematics \\ % \hfill (Received 00 00 2010)\\
Ibn Tofail University   \\ %\hfill (Revised  00 00 2010)\\
Kenitra,Morrocco\\
}
\keywords{special functions, contiguous hypergeometric functions}

% For each author, add one of the following

\date{\today}
  
\begin{abstract}
The aim of this paper is to give, using some contiguous relations, the asymptotic behaviour of some linear combination of two symmetric contiguous hypergeometric functions, under some conditions of their parameters.

\end{abstract}

\maketitle

\section{Introduction, notations and statement of the main results }

Hypergeometric function  belongs to an important class of special functions. They are very useful in many fields of physical and astronomical researches, in mathematical analysis and its application. They have significant properties: identities of special values, and transformation formulas. These properties have many intereting applications in combinatorial analsis and arithmetic geometric.

The hypergeometric series is defined by: $${}_2F_1\left( a,b; c; z \right) = \Sigma_{n \geq 0} \frac{(a)_n (b)_n}{(c)_n} \frac{z^n}{n!}
$$
with $(a)_n:=a(a+1) \cdots (a+n-1)$ for $n \geq 1$,  $ (a)_{0} = 1  $

where it is assumed that $ c \neq 0,-1,-2,-3,...  $ so that no zero factors appear in the denominators of the terms of the series. The variable is z, and a,b and c are called the parameters of the function. 

  The first systematic and thorough study of hypergeometric series was carried out by Gauss in his paper of $1820  $. Many people: Ernst Kummer, Jacobi, Riemann, E.W Barnes, RH Mellin gived different approches and interesting results about hypergeometric functions. 

An hypergeometric function is called contiguous to the other $F(a,b,c,z)$ if one, two, or three of the parameters $ a,b,c $ are increased or diminished by an integer.

In his $1820$ paper, Gauss first introduced the terminology "contiguous function". He definened two hypergeometric function to be contiguous if they have the same power series variable and if two of the parameters are pairwise equal and the third pair differs by $- 1 $. He found that every three contiguous hypergeometric functions are linearly related. Then he listed $15$ such linear relations (which are now called "contiguous relations"), and he used them to derive basic formulas.
 
These contiguous hypergeometric functions have some interesting applications.\\
Indeed, contiguous relations gives an intetwining correspondance between Lie algebras and special functions.They are also very useful in the derivation of summation and transformation formulas for hypergeometric series. 

Recently, many people studies contiguous hypergeometric functions and 
gives some new properties, consequences and interesting applications of them. See for instance [5],[6],[7],[8],[9],[10],[11].

in [1],the autors used well-known forms of the Gauss function to prove some simple identities relating shifted hypergeometric functions together with their derivatives. Then they combined these identities to obtain more general new identities.

Usually, while solving many interesting problems related to mathematic analysis, we are led to establish a precise asymptotic behavior of some functions. In this paper, using some contiguous relations of [1] , we are able to give, under some conditions of their parameters, the asymptotic behaviour of linear combination of two symmetric contiguous hypergeometric functions.

Now we may state the main theorem of this paper:\\

\begin{thm} 

For $(a-b), \alpha, \beta \in \mathbb{N} $ and $ a,b,c \in \mathbb{C} $ with:\\
\begin{equation*} Re(a+b+\alpha+\beta-c-1)> 0 \end{equation*} we have\\
\begin{equation*}
{lim}_{Z\rightarrow1}(1-Z)^{-c+a+b+\alpha+\beta-1}( {(a)}_{\alpha}{(b)}_{\beta}{F}_{21}(a+\alpha,b+\beta;c;z)\\
-{(a)}_{\beta}{(b)}_{\alpha}{F}_{21}(a+\beta,b+\alpha;c;z)) \end{equation*}
\begin{equation*}
=\frac{\Gamma(c)}{\Gamma(a)\Gamma(b)}\Gamma(a+b+\alpha+\beta-c-1)(a-b)(\alpha-\beta)
\end{equation*}
 \end{thm}

\section{Proof of the main theorem }

\begin{proof}
Denote by $ {F}_{\alpha,\beta,\gamma}(Z) $ the hypergeometric function:\\\

 $ {F}_{21}(a+\alpha,b+\beta;c+\gamma;z) $ for $ \alpha, \beta, \gamma \in \mathbb{N} $ and $ a,b,c \in \mathbb{C} $\\

Then we have the following identitie: (see for instance [1] p:407)

\begin{equation*}
{F}_{\alpha,\beta,\gamma}(Z)= {F}_{\alpha-1,\beta,\gamma}(Z)+ \frac{b+\beta}{c+\gamma}Z{F}_{\alpha,\beta+1,\gamma+1}(Z)
\end{equation*}

Then we can rwite:
\begin{equation*}
{F}_{\alpha,\beta,0}(Z)= {F}_{\alpha-1,\beta,0}(Z)+ \frac{b+\beta}{c}Z{F}_{\alpha,\beta+1,1}(Z)
\end{equation*}

And similarly:

\begin{equation*}
{F}_{\beta,\alpha,0}(Z)= {F}_{\beta-1,\alpha,0}(Z)+ \frac{b+\alpha}{c}Z{F}_{\beta,\alpha+1,1}(Z)
\end{equation*}

\begin{equation*}
{(a)}_{\alpha}{(b)}_{\beta}{F}_{\alpha,\beta,0}(Z)= {(a)}_{\alpha}{(b)}_{\beta}{F}_{\alpha-1,\beta,0}(Z)+ {(a)}_{\alpha}{(b)}_{\beta+1}\frac{Z}{c}{F}_{\alpha,\beta+1,1}(Z)
\end{equation*}

And similarly:

\begin{equation*}
{(a)}_{\beta}{(b)}_{\alpha}{F}_{\beta,\alpha,0}(Z)= {(a)}_{\beta}{(b)}_{\alpha}{F}_{\beta-1,\alpha,0}(Z)+ {(a)}_{\beta}{(b)}_{\alpha+1}\frac{Z}{c}{F}_{\beta,\alpha+1,1}(Z)
\end{equation*}

And

\begin{equation*}
{(a)}_{\alpha}{(b)}_{\beta}{F}_{\alpha,\beta,0}(Z)= {(a)}_{\alpha}{(b)}_{\beta}{F}_{\alpha-1,\beta,0}(Z)+ 
\end{equation*}
\begin{equation*}
{(a)}_{\alpha}{(b)}_{\beta+1}\frac{Z}{c}[{F}_{\alpha-1,\beta+1,1}(Z)+ \frac{b+\beta+1}{c+1}Z {F}_{\alpha,\beta+2,2} ]
\end{equation*}

And similarly:

\begin{equation*}
{(a)}_{\beta}{(b)}_{\alpha}{F}_{\beta,\alpha,0}(Z)= {(a)}_{\beta}{(b)}_{\alpha}{F}_{\beta-1,\alpha,0}(Z)+\end{equation*}
\begin{equation*} {(a)}_{\beta}{(b)}_{\alpha+1}\frac{Z}{c}[{F}_{\beta-1,\alpha+1,1}(Z)+ \frac{b+\alpha+1}{c+1}Z {F}_{\beta,\alpha+2,2}]
\end{equation*}

We get:

\begin{equation*}
{(a)}_{\alpha}{(b)}_{\beta}{F}_{\alpha,\beta,0}(Z)= {(a)}_{\alpha}{(b)}_{\beta}{F}_{\alpha-1,\beta,0}(Z)+
\end{equation*}
\begin{equation*}
\frac{Z}{c} {(a)}_{\alpha}{(b)}_{\beta+1}{F}_{\alpha-1,\beta+1,1}(Z)+ \frac{Z^{2}}{c(c+1)}{(a)}_{\alpha}{(b)}_{\beta+2} {F}_{\alpha,\beta+2,2}(Z) 
\end{equation*}

And similarly:

\begin{equation*}
{(a)}_{\beta}{(b)}_{\alpha}{F}_{\beta,\alpha,0}(Z)= {(a)}_{\beta}{(b)}_{\alpha}{F}_{\beta-1,\alpha,0}(Z)+
\end{equation*}
\begin{equation*}
\frac{Z}{c} {(a)}_{\beta}{(b)}_{\alpha+1}{F}_{\beta-1,\alpha+1,1}(Z)+ \frac{Z^{2}}{c(c+1)} {(a)}_{\beta}{(b)}_{\alpha+2} {F}_{\beta,\alpha+2,2}(Z)
\end{equation*}

But we also have:

\begin{equation*}
{F}_{\alpha,\beta+2,2}(Z)= {F}_{\alpha-1,\beta+2,2}(Z)+ \frac{b+\beta+2}{c+2}Z{F}_{\alpha,\beta+3,3}(Z)
\end{equation*}

And similarly:

\begin{equation*}
{F}_{\beta,\alpha+2,0}(Z)= {F}_{\beta-1,\alpha+2,2}(Z)+ \frac{b+\alpha+2}{c+2}Z{F}_{\beta,\alpha+3,3}(Z)
\end{equation*}

Wich gives:

\begin{equation*}
{(a)}_{\alpha}{(b)}_{\beta}{F}_{\alpha,\beta,0}(Z)= {(a)}_{\alpha}{(b)}_{\beta}{F}_{\alpha-1,\beta,0}(Z)+\frac{Z}{c} {(a)}_{\alpha}{(b)}_{\beta+1}{F}_{\alpha-1,\beta+1,1}(Z)+ 
\end{equation*}
\begin{equation*}
\frac{Z^{2}}{c(c+1)}{(a)}_{\alpha}{(b)}_{\beta+2} {F}_{\alpha-1,\beta+2,2}(Z) 
\end{equation*}
\begin{equation*}
+ \frac{Z^{3}}{c(c+1)(c+2)}{(a)}_{\alpha}{(b)}_{\beta+3} {F}_{\alpha,\beta+3,3}(Z) 
\end{equation*}

And similarly:

\begin{equation*}
{(a)}_{\beta}{(b)}_{\alpha}{F}_{\beta,\alpha,0}(Z)= {(a)}_{\beta}{(b)}_{\alpha}{F}_{\beta-1,\alpha,0}(Z)+\frac{Z}{c} {(a)}_{\beta}{(b)}_{\alpha+1}{F}_{\beta-1,\alpha+1,1}(Z)+ 
\end{equation*}
\begin{equation*}
\frac{Z^{2}}{c(c+1)} {(a)}_{\beta}{(b)}_{\alpha+2} {F}_{\beta-1,\alpha+2,2}(Z)
\end{equation*}
\begin{equation*}
+ \frac{Z^{3}}{c(c+1)(c+2)}{(a)}_{\beta}{(b)}_{\alpha+3} {F}_{\beta,\alpha+3,3}(Z) 
\end{equation*}

And so on

By iterating the process $ (a-b) $ times we finally get:
\begin{equation*}
{(a)}_{\alpha}{(b)}_{\beta}{F}_{\alpha,\beta,0}(Z)= \sum_{i=0}^{a-b-1}\frac{Z^{i}}{(c)_{i}}{(a)}_{\alpha}{(b)}_{\beta+i}{F}_{\alpha-1,\beta+i,i}(Z)+
\end{equation*}
\begin{equation*}
\frac{Z^{a-b}}{(c)_{a-b}} {(a)}_{\alpha}{(b)}_{\beta+a-b}{F}_{\alpha,\beta+a-b,a-b}(Z)
\end{equation*}

And similarly:

\begin{equation*}
{(a)}_{\beta}{(b)}_{\alpha}{F}_{\beta,\alpha,0}(Z)= \sum_{i=0}^{a-b-1}\frac{Z^{i}}{(c)_{i}}{(a)}_{\beta}{(b)}_{\alpha+i}{F}_{\beta-1,\alpha+i,i}(Z)+
\end{equation*}
\begin{equation*}
\frac{Z^{a-b}}{(c)_{a-b}} {(a)}_{\beta}{(b)}_{\alpha+a-b}{F}_{\beta,\alpha+a-b,a-b}(Z)
\end{equation*}

But we also have:
\begin{equation*}
{(a)}_{\alpha}{(b)}_{\beta+a-b}{F}_{\alpha,\beta+a-b,a-b}(Z)=  {(a)}_{\beta}{(b)}_{\alpha+a-b}{F}_{\beta,\alpha+a-b,a-b}(Z)
\end{equation*}
\begin{equation*}
= \frac{\Gamma(a+\alpha) \Gamma(a+\beta)}{\Gamma(a) \Gamma(b)}{F}(a+\alpha,a+\beta,a-b,Z)
\end{equation*}

And for each $ i $  in $ [0,a-b-1]$ we have:

\begin{equation*}
\frac{Z^{i}}{(c)_{i}}{(a)}_{\alpha}{(b)}_{\beta+i}{F}_{\alpha-1,\beta+i,i}(Z) \sim_{1}
\end{equation*}
\begin{equation*}
 (1-Z)^{c-a-b-\alpha-\beta+1}(a+\alpha-1)\frac{\Gamma(c)}{\Gamma(a) \Gamma(b)} \Gamma(a+b+\alpha+\beta-c-1)
\end{equation*}

And similarly

\begin{equation*}
\frac{Z^{i}}{(c)_{i}}{(a)}_{\beta}{(b)}_{\alpha+i}{F}_{\beta-1,\alpha+i,i}(Z)\sim_{1}
\end{equation*}
\begin{equation*}
 (1-Z)^{c-a-b-\alpha-\beta+1}(a+\beta-1)\frac{\Gamma(c)}{\Gamma(a) \Gamma(b)} \Gamma(a+b+\alpha+\beta-c-1)
\end{equation*}
Wich gives the desired result;
\end{proof}

\bibliographystyle{amsplain}

\end{document}